\documentclass[12pt]{amsart}
\usepackage{amscd}
\usepackage{amssymb}

\def\frk{\frak}               

\def\pp{{\frk p}}

\def\mm{{\frk m}}

\def\Phi{{\frk n}}
\def\Phi{{\frk N}}
\def\opn#1#2{\def#1{\operatorname{#2}}} 
\opn\chara{char} \opn\length{\ell} \opn\pd{pd} \opn\rk{\lk}\opn\link{link}
\opn\projdim{proj\,dim} \opn\injdim{inj\,dim} \opn\rank{rank}
\opn\depth{depth} \opn\and{and} \opn\grade{grade}
\opn\height{height} \opn\embdim{emb\,dim} \opn\codim{codal}

\opn\Tr{Tr} \opn\bigrank{big\,rank}
\opn\superheight{superheight}\opn\lcm{lcm}
\opn\trdeg{tr\,deg}%
\opn\reg{reg} \opn\lreg{lreg} \opn\ini{in}
\opn\div{div} \opn\Div{Div} \opn\cl{cl} \opn\Cl{Cl}
\opn\Spec{Spec} \opn\Supp{Supp} \opn\supp{supp} \opn\Sing{Sing}
\opn\Ass{Ass} \opn\Min{Min}
\opn\Ann{Ann} \opn\Rad{Rad} \opn\Soc{Soc}
\opn\Im{Im}
 \opn\Ker{Ker} \opn\Coker{Coker} \opn\Am{Am} \opn \inf{inf}
\opn\Hom{Hom} \opn\Tor{Tor} \opn\Ext{Ext} \opn\End{End} \opn\cd{cd}\opn\Var{Var}
\opn\Aut{Aut} \opn\id{id}

\opn\nat{nat}
\opn\pff{pf}
\opn\Pf{Pf} \opn\GL{GL} \opn\SL{SL} \opn\mod{mod} \opn\ord{ord}
\opn\cl{cl} \opn\conv{conv} \opn\ext{ext} \opn\rad{rad}
\opn\star{star} \opn\red{red}\opn\H{H}
\opn\aff{aff} \opn\con{conv} \opn\relint{relint} \opn\st{st}
\opn\lk{lk} \opn\cn{cn} \opn\core{core} \opn\vol{vol}
\opn\link{link} \opn\star{star}
\opn\gr{gr}

\def\pot#1#2{#1[\kern-0.28ex[#2]\kern-0.28ex]}

\opn\dirlim{\underrightarrow{\lim}}
\opn\inivlim{\underleftarrow{\lim}}
\let\dirsum=\oplus
\let\tensor=\otimes
\let\iso=\cong

\let\to=\rightarrow

\def\Implies{\ifmmode\Longrightarrow \else
     \unskip${}\Longrightarrow{}$\ignorespaces\fi}
\def\implies{\ifmmode\Rightarrow \else
     \unskip${}\Rightarrow{}$\ignorespaces\fi}
\def\iff{\ifmmode\Longleftrightarrow \else
     \unskip${}\Longleftrightarrow{}$\ignorespaces\fi}

\let\:=\colon
\newtheorem{Theorem}{Theorem}[section]
\newtheorem{Lemma}[Theorem]{Lemma}

\newtheorem{Proposition}[Theorem]{Proposition}
\newtheorem{Remark}[Theorem]{Remark}

\newtheorem{Example}[Theorem]{Example}

\newtheorem{Definition}[Theorem]{Definition}

\newtheorem{Fact}[Theorem]{Fact}
\let\epsilon\varepsilon
\let\phi=\varphi
\let\kappa=\varkappa
\def\qed{\ifhmode\textqed\fi
   \ifmmode\ifinner\quad\qedsymbol\else\dispqed\fi\fi}
\def\textqed{\unskip\nobreak\penalty50
    \hskip2em\hbox{}\nobreak\hfil\qedsymbol
    \parfillskip=0pt \finalhyphendemerits=0}
\def\dispqed{\rlap{\qquad\qedsymbol}}

\opn\dis{dis}
\def\pnt{{\raise0.5mm\hbox{\large\bf.}}}

\begin{document}
\title{Sequentially Cohen--Macaulayness of bigraded modules}

\author{  Ahad Rahimi}

\subjclass[2000]{ 16W50, 13C14, 13D45, 16W70.
The author was in
part supported by a grant from IPM (No. 900130058)}

 \address{ Ahad Rahimi, Department of Mathematics, Faculty of Science, Razi University, Baghe abrisham, Kermanshah,
 Iran and
 School of Mathematics, Institute for Research in Fundamental Sciences
(IPM), P. O. Box: 19395-5746, Tehran, Iran.
}\email{ahad.rahimi@razi.ac.ir}

\begin{abstract}
Let $K$ be a field, $S=K[x_1,\ldots,x_m, y_1,\ldots,y_n]$ be a standard bigraded polynomial ring and $M$ a finitely generated bigraded $S$-module. In this paper we study sequentially Cohen--Macaulayness of $M$ with respect to $Q=(y_1,\ldots,y_n)$. We characterize the sequentially Cohen--Macaulayness of $L\tensor_KN$
with respect to $Q$ as an $S$-module when $L$ and $N$ are non-zero finitely generated
graded modules over $K[x_1, \dots, x_m]$ and $K[y_1, \dots, y_n]$, respectively.
 All hypersurface rings that are sequentially Cohen--Macaulay with
respect to $Q$  are classified.
\end{abstract}

\maketitle
\section*{Introduction}
In \cite {St}  Stanley introduced the notion of
 sequentially Cohen--Macaulayness for graded modules. This concept has been then studied by several authors, we refer to  \cite{CN}, \cite{F}, \cite{HS}, \cite{PR}, \cite{NR}, \cite{S}, \cite{TY}.
  In this paper we define sequentially Cohen-Macaulayness
for bigraded modules, and introduce some new algebraic invariants which it makes sense to study in this case.
 We let $S=K[x_1, \dots, x_m, y_1, \dots, y_n]$ be a standard bigraded
 polynomial ring over a field $K$,  $M$ a finitely
 generated bigraded $S$-module.  We set  $Q=(y_1, \dots,  y_n)$.    In \cite{AR2}, $M$ is called to be Cohen--Macaulay with respect to $Q$
 if  $\grade(Q, M)=\cd(Q, M)$,
 where  $\cd(Q, M)$ denotes the cohomological dimension of $M$ with respect to $Q$.

We call a finite filtration $\mathcal{F}$: $0=M_0\varsubsetneq M_1
\varsubsetneq
 \dots  \varsubsetneq M_r=M$
 of $M$ by
bigraded submodules $M$ a Cohen--Macaulay filtration with
respect to $Q$ if
 \begin{itemize}
\item[{(a)}] Each quotient $M_i/M_{i-1}$ is  Cohen--Macaulay with respect to $Q$;
 \item[{(b)}] $0 \leq \cd(Q, M_1/M_0)<\cd(Q, M_2/M_1)< \dots< \cd(Q, M_r/M_{r-1})$.
 \end{itemize}
 If $M$ admits a Cohen--Macaulay filtration with respect to $Q$, then we say that $M$ is sequentially Cohen--Macaulay with respect to $Q$.
   The usual notion of sequentially Cohen-Macaulayness
 results from our definition if
we assume $P=0$.

A finite filtration $
\mathcal{D}:  0=D_0\varsubsetneq D_1
\varsubsetneq
 \dots  \varsubsetneq  D_r=M
 $
 of $M$ by bigraded submodules, is called the
dimension filtration of $M$ with respect to $Q$ if $D_{i-1}$ is the
largest bigraded submodule of $D_i$ for which $\cd(Q,D_{i-1})<\cd(Q, D_i)$,
for all $i=1, \dots, r$.
In the preliminary section we show that if $M$ is sequentially Cohen--Macaulay with respect to $Q$, then the filtration $\mathcal{F}$ is uniquely determined and it is just the dimension filtration of $M$ with respect to $Q$, that is, $\mathcal{F}=\mathcal{D}$. We explicitly describe the structure of the submodules $D_i$ in \cite{PR}. In the same section, it is also shown that if $M$ is sequentially
Cohen--Macaulay with respect to $Q$ with $\grade(Q, M)>0$  and  $|K|=\infty$, then there exists a
bihomogeneous $M$-regular element $y\in Q$ of degree $(0,1)$ such
that $M/yM$ is sequentially Cohen--Macaulay with respect to $Q$, too. An example is given to show that the converse does not hold in general.

Let $K[x]=K[x_1, \dots, x_m]$ and $K[y]=K[y_1, \dots, y_n]$.  In Section 2,  we consider $L\tensor_KN$ as $S$-module where $L$ and
$N$ are two non-zero finitely generated graded modules over $K[x]$
and $K[y]$, respectively. We characterize the sequentially Cohen--Macaulayness of $L\tensor_KN$
with respect to $Q$ as follows: $L\tensor_KN$ is a sequentially Cohen--Macaulay
with respect to $Q$ if and only if $N$ is sequentially Cohen--Macaulay $K[y]$-module.

In  the last section, we let $f\in S$ be a bihomogeneous element of
degree $(a, b)$ and consider the hypersurface ring $R=S/fS$. Notice
that if $a, b>0$,  we have $\grade(Q, R)=n-1$ and $\cd(Q, R)=n$;  hence $R$ is not
Cohen--Macaulay with respect to $Q$. Thus,  it is natural to
ask whether $R$ is sequentially Cohen--Macaulay with respect to $Q$.
 We classify all hypersurface rings that are sequentially
Cohen--Macaulay with respect to $Q$. In fact, we show:
 $R$ is sequentially Cohen--Macaulay
 with respect to $Q$ if and only
if $f=h_1h_2$ where $\deg(h_1)=(a, 0)$ with $a\geq 0$  and
$\deg(h_2)=(0, b)$ with $b\geq 0$.

\section{Preliminaries}
 Let $K$ be a field and $S=K[x_1, \dots, x_m, y_1, \dots, y_n]$ be a standard
bigraded
 polynomial ring over $K$. In other words, $\deg x_i=(1,0)$ and $\deg y_j=(0, 1)$  for all $i$ and $j$.
 We set
$P=(x_1, \dots, x_m)$ and $Q=(y_1, \dots,  y_n)$.
Let $M$ be a
finitely
 generated bigraded $S$-module.
We denote by $\cd(Q, M)$ the {\em cohomological dimension of $M$ with
respect to $Q$ } which is the largest integer $i$ for which $H^i_ Q
(M)\neq 0$. Notice that $0\leq \cd(Q, M) \leq n.$

\begin{Definition}
\label{cd}
{\em   We say $M$ is {\em Cohen--Macaulay with respect to $Q$} if we have only one non vanishing local cohomology module with respect to $Q$. This was previously called under a different name in \cite{AR2}. Namely, {\em relative Cohen--Macaulay with respect to $Q$},  which for simplicity's sake, we omit the word "relative". }
\end{Definition}
We recall the following facts which will be used in the sequel.
\begin{Fact}{\em
\label{cd}
Let $M$ be a finitely generated bigraded $S$-module. Then
 \begin{itemize}
 \item[{(a)}]
 $\cd(P, M)=\dim M/QM$  and $\cd(Q, M)=\dim M/PM$,  see \cite[Formula 3]{AR2}.
\item[{(b)}] $\grade(Q, M)\leq \dim M-\cd(P, M)$, and the equality holds if $M$ is Cohen--Macaulay, see \cite[Formula 5]{AR2}.
\item[{(c)}] The exact sequence  $ 0
\rightarrow M' \rightarrow M \rightarrow M'' \rightarrow 0$ of finitely generated bigraded $S$-modules yields $ \cd(Q,M)=\max\{\cd(Q, M'), \cd(Q,M'')\}$,
see the general version of \cite[Proposition 4.4]{CJR}.
\item[{(d)}]   $\cd(Q,M)=\max \{\cd (Q, S/{\pp}): \pp \in \Ass(M)\}$,
see the general version of \cite[Corollary 4.6]{CJR}.
\end{itemize}}
\end{Fact}




\begin{Definition}
\label{seq1}{\em We call a finite filtration $\mathcal{F}$:
$0=M_0\varsubsetneq M_1 \varsubsetneq
 \dots  \varsubsetneq M_r=M$
 of $M$ by
bigraded submodules a {\em Cohen--Macaulay filtration with
respect to $Q$} if
 \begin{itemize}
\item[{(a)}] Each quotient $M_i/M_{i-1}$ is Cohen--Macaulay with respect to $Q$;
 \item[{(b)}] $0 \leq \cd(Q, M_1/M_0)<\cd(Q, M_2/M_1)< \dots< \cd(Q, M_r/M_{r-1})$.
 \end{itemize}
If $M$ admits a Cohen--Macaulay filtration with respect to $Q$, then we say $M$ is
{\em sequentially Cohen--Macaulay with respect to $Q$}.  }
\end{Definition}
Observe that the ordinary definition of sequentially Cohen--Macaulay modules results from our definition if we assume $P=0$.
\begin{Remark}
\label{dimeq}{\em
By applying Fact \ref{cd}(a) to the exact sequences $0\to M_{i-1}\to M_i \to M_i/M_{i-1} \to 0$ one immediately has
  $\cd(Q, M_i)=\cd(Q, M_i/M_{i-1})$ for $i=1, \dots, r.$}
\end{Remark}

\begin{Example}{\em
 Cohen--Macaulay modules with respect to $Q$ are obvious
examples of sequentially Cohen--Macaulay modules with respect to $Q$. Any module $M$ such that $\cd(Q, M)\leq 1$ is sequentially Cohen--Macaulay with
respect to $Q$.  To see this, we may assume that $M$ is not Cohen--Macaulay with respect to $Q$. Thus  $\grade(Q,M)=0$ and $\cd(Q,
M)=1$. The filtration $0=M_0 \varsubsetneq M_1 \varsubsetneq
M_2=M$ where $M_1=H^0_Q(M)$ is a Cohen-Macaulay filtration
with respect to $Q$. }
\end{Example}

In the following, we show that the filtration
 $\mathcal{F}$  given
in Definition \ref{seq1}   is unique. To do so we
need some preparation.
\begin{Lemma}
\label{uni}
There is a unique largest bigraded submodule $N$ of $M$ for which $\cd(Q,N)<\cd(Q,
M)$.
\end{Lemma}
\begin{proof}
Let  $\sum$ be the set of all bigraded submodules $L$ of $M$ such that
 $\cd(Q,L)<\cd(Q, M)$. As $M$ is a Noetherian $S$-module, $\sum$
has a maximal element with respect to inclusion, say $N$. Let $T$ be
an arbitrary element in $\sum$. Fact \ref{cd}(c) implies $\cd(Q,
T+N)<\cd(Q, M)$,  hence the maximality of $N$ yields $T\subseteq N$.
\end{proof}
\begin{Definition}
\label{1} {\em A filtration $\mathcal{D}$: $0=D_0\varsubsetneq D_1
\varsubsetneq
 \dots  \varsubsetneq  D_r=M$ of $M$ by bigraded submodules  is called the
 { \em dimension filtration of $M$ with respect to $Q$ } if $D_{i-1}$ is the
largest bigraded submodule of $D_i$ for which $\cd(Q, D_{i-1})<\cd(Q, D_i)$
for all $i=1, \dots, r$.}
\end{Definition}

The dimension filtration introduced by Schenzel \cite{S}  is thus a
dimension filtration with respect to the maximal ideal $\mm=P+Q$. A filtration $\mathcal{D}$ as in Definition \ref{1} is unique by Lemma \ref{uni}.
In order to prove the  uniqueness of an $\mathcal{F}$ as in Definition \ref{seq1}, we
will show that  $\mathcal{F}=\mathcal{D}$. In \cite{JR},   $M$ is called to be
 {\em relatively unmixed}  with respect to $Q$ if $\cd(Q, M)=\cd(Q, S/\pp)$
for all $\pp \in \Ass(M).$
\begin{Lemma}
 \label{CM}
 Let $N$  be a non-zero
 bigraded submodule of $M$. If $M$ is
  Cohen--Macaulay with respect to $Q$, then $\cd(Q, N)=\cd(Q, M)$.
\end{Lemma}
\begin{proof}
Since $M$ is Cohen--Macaulay with respect to $Q$, it
follows from \cite[Corollary 1.11]{JR} that $M$ is relatively
unmixed with respect to $Q$, i.e., $\cd(Q, M)=\cd(Q, S/\pp)$ for
all $\pp \in \Ass(M)$. As $N$ is a non-zero submodule of $M$ we have
$\Ass(N)\neq\emptyset$  and  $\Ass(N)\subseteq \Ass(M)$. Thus,  Fact \ref{cd}(d) implies
\[
\cd(Q,N)  =  \max \{\cd (Q, S/{\pp}): \pp \in
\Ass(N)\}=\cd(Q, M),
\]
 as desired.
\end{proof}
\begin{Proposition}
\label{unique} Let $\mathcal{F}$ be a Cohen--Macaulay
filtration of $M$
 with respect to $Q$ and
 $\mathcal{D}$ be the dimension filtration of $M$ with respect to $Q$.
 Then
$\mathcal{F}=\mathcal{D}$.
\end{Proposition}
\begin{proof}
Let $\mathcal{F}$: $0=M_0\varsubsetneq M_1 \varsubsetneq  \dots
\varsubsetneq M_r=M$ and  $\mathcal{D}$: $0=D_0\varsubsetneq D_1
\varsubsetneq \dots \varsubsetneq D_s=M$; we will show that $r=s$
and $M_i=D_i$ for all $i$. By Remark \ref{dimeq}
 we have $\cd(Q, M_{i-1})<\cd(Q, M_{i})$ for all $i=1, \dots, r.$ Hence,
  Definition \ref{1} says that  $M_{r-1}\subseteq D_{s-1}$. Assume
$M_{r-1}\varsubsetneq D_{s-1}$. Thus $D_{s-1}/M_{r-1}$ is a non-zero
submodule of $M/M_{r-1}$. Since $M/M_{r-1}$ is
Cohen--Macaulay with respect to $Q$, it follows from Lemma \ref{CM}
that  $\cd(Q, D_{s-1}/M_{r-1})=\cd(Q,M/M_{r-1})=\cd(Q, M)$, where the second equality is yielded by
 Remark \ref{dimeq}. Now applying Fact \ref{cd}(c) to the exact sequence $0 \rightarrow M_{r-1} \rightarrow D_{s-1}
\rightarrow D_{s-1}/M_{r-1} \rightarrow 0$  yields $\cd(Q,
D_{s-1})=\cd(Q, M)$, a contradiction. Thus $M_{r-1}=D_{s-1}$.
Continuing in this way,  we get $r=s$ and $M_i=D_i$ for all $i$.
Therefore, $\mathcal{F}=\mathcal{D}$.
\end{proof}

 We end this section with the following proposition which gives us a class of sequentially Cohen--Macaulay with respect to $Q$. First we have the following

\begin{Lemma}
\label{unmixed} Let $M$ be sequentially Cohen--Macaulay with respect to $Q$. If $M$ is
relatively unmixed  with respect to $Q$, then $M$ is Cohen--Macaulay with respect to $Q$.
\end{Lemma}
\begin{proof}
Let $0=M_0\varsubsetneq M_1 \varsubsetneq  \dots  \varsubsetneq
M_r=M$ be the Cohen--Macaulay filtration with respect to
$Q$.
 By Fact \ref{grade}  we have $\grade(Q, M)=\grade(Q, M_1)$. Since $M_1$
 is Cohen--Macaulay with respect to $Q$, it follows  from \cite[Corollary 1.11]{JR} that
  $M_1$ is relatively unmixed  with respect to $Q$.
   Thus  $\grade(Q, M)=\grade(Q, M_1)=\cd(Q, M_1)=\cd(Q, S/\pp)$ for all $\pp \in \Ass(M_1)$.
 As  $M$ is relatively unmixed  with respect to $Q$ and $\Ass(M_1)\subseteq \Ass(M)$, we have
  $\grade(Q, M)=\cd(Q, M)$, as desired.
\end{proof}

\begin{Proposition}
\label{regular} Suppose $\grade(Q, M)>0$  and  $|K|=\infty$. If $M$ is sequentially
Cohen--Macaulay with respect to $Q$, then there exists a
bihomogeneous $M$-regular element $y\in Q$ of degree $(0,1)$ such
that $M/yM$ is sequentially Cohen--Macaulay with respect to $Q$.
\end{Proposition}
\begin{proof}
We assume that $M$ is sequentially Cohen--Macaulay with
respect to $Q$ and let $\mathcal{F}$: $0=M_0\varsubsetneq M_1
\varsubsetneq \dots \varsubsetneq M_r=M$ be the
Cohen--Macaulay filtration with respect to $Q$. Since $\grade(Q, M)=\grade(Q, M_1)=\cd(Q, M_1)>0$, it follows that $\grade(Q, M_i/M_{i-1})=\cd(Q, M_i/M_{i-1})>0$ for all $i$. We set $N_i=M_i/M_{i-1}$.
Thus by \cite[Corollary, 3.5]{AR2} which is also valid for finitely many modules which are Cohen--Macaulay with respect to $Q$ and have positive cohomological dimension with respect to $Q$,
there exists a bihomogeneous element $y\in Q$ of degree $(0, 1)$ such that $y$ is $N_i$-regular for all $i$ and  $\overline{N_i}$ is Cohen--Macaulay with respect to $Q$ with $\cd(Q, \overline{N_i})=\cd(Q, N_i)-1$. Here $\overline{L}=L/yL$ for any $S$-module $L$.

   Consider the exact sequence $0 \rightarrow M_{i-1} \rightarrow M_i \rightarrow N_i \rightarrow 0$ for all $i$. Since $y$ is regular on $N_i$ for all $i$, it follows that $\Tor^S_1(S/yS, N_i)=0$ for all $i$. Hence, we get the following exact sequence $0 \rightarrow \overline{M_{i-1}} \rightarrow \overline{M_{i}} \rightarrow \overline{N_i} \rightarrow 0$ for all $i$.
Now the filtration $\mathcal{G}$: $0=\overline{M_0} \varsubsetneq \overline{M_1}
\varsubsetneq \dots \varsubsetneq \overline{M_r}=M/yM$ is the Cohen--Macaulay filtration for $M/yM$ with respect to $Q$. In fact,   $\overline{M_i} /\overline{M_{i-1}}\iso \overline{N_i}$ and $\grade(Q, \overline{N_i}) = \grade(Q, N_i)-1= \cd(Q,  N_i)-1= \cd(Q, \overline{N_i})$. Hence,  $\grade(Q, \overline{M_i} /\overline{M_{i-1}}) = \cd(Q, \overline{M_i} /\overline{M_{i-1}})$. As  $\cd(Q, M_i/M_{i-1})<\cd(Q, M_{i+1}/M_{i})$  for all $i$, we have $\cd(Q,\overline{M_i} /\overline{M_{i-1}})<\cd(Q, \overline{M_{i+1}}/\overline{M_{i}})$ for all $i$.
 \end{proof}
The following example shows that the converse of Proposition \ref{regular} does not hold in general.

\begin{Example}{\em
Consider the hypersurface ring $R=K[x_1, x_2, y_1, y_2]/(f)$ where $f=x_1y_1+x_2y_2.$ One has $\grade(Q, R)=1$ and $\cd(Q, R)=2.$ By \cite[Lemma 3.4]{AR2} there exists a bihomogeneous $R$-regular element $y\in Q$ of degree $(0, 1)$ such that $\cd(Q, R/yR)=\cd(Q, R)-1=1$ and of course $\grade(Q, R/yR)=\grade(Q, R)-1=0$. Hence $R/yR$ is sequentially Cohen--Macaulay with respect to $Q$. On the other hand,  $R$ is not sequentially Cohen--Macaulay with respect to $Q$. Indeed, $\Ass(R)=\{(f)\}$ and  $\cd(Q, R)=\cd(Q, S/(f))$ says that $R$ is relatively unmixed with respect to $Q$. If $R$ is sequentially Cohen--Macaulay with respect to $Q$, then by Lemma \ref{unmixed} $R$ is Cohen--Macaulay with respect to $Q$, a contradiction.
}
\end{Example}

\section{Sequentially Cohen--Macaulayness of $L\tensor_KN$ with respect to $Q$}
In this section,  we characterize the sequentially Cohen--Macaulayness of
$L\tensor_KN$ with respect to $Q$ as $S$-module where $L$ and $N$ are two non-zero finitely generated
graded modules over $K[x]$ and $K[y]$, respectively. For the bigraded $S$-module $M$ we define the bigraded Matlis-dual of $M$ to be $M^\vee$ where the $(-i, -j)$th bigraded components of $M^\vee$ is given by $\Hom_K(M_{(i, j)}, K).$ We set $M_k = M_{(k,*)} = \dirsum_jM_{(k,j)}$
and consider it as a finitely generated graded $K[y]$-module.

\begin{Lemma}
\label{cmd}
Let $M$ be a finitely generated bigraded $S$-module. If $M$ is Cohen--Macaulay with respect to $Q$ with $\cd(Q, M)=q$, then $\big( H^{q}_Q(M)^\vee\big)_{(k, *)}$ is a finitely generated Cohen--Macaulay $K[y]$-module of dimension $q$ for all $k$.
\end{Lemma}
\begin{proof} Note that
\begin{eqnarray*}
\big( H^{i}_Q(M)^\vee\big)_{(k, *)} &\iso &  \big( H^{i}_Q(M)_{(-k, *)}\big)^\vee\\
                                    &\iso &    \big( H^{i}_{(y_1, \dots, y_n)}(M_{(-k, *)})\big)^\vee\\
                                    &\iso &  \Ext^{n-i}_{K[y]}(M_{(-k, *)}, K[y](-n)).
\end{eqnarray*}

Since $M$ is Cohen--Macaulay with respect to $Q$ with $\cd(Q, M)=q$, it follows from \cite[Proposition 1.2]{AR2} that $M_{(-k, *)}$ is a Cohen--Macaulay $K[y]$-module of dimension $q$, and the conclusion follows immediately.
\end{proof}
\begin{Lemma}
\label{reg}
 Let $M$ be sequentially Cohen--Macaulay with respect to $Q$
 with Cohen--Macaulay filtration $\mathcal{F}$:
 $0=M_0\varsubsetneq M_1 \varsubsetneq  \dots \varsubsetneq M_r=M$ with respect to $Q$.
 Then we have
$H^{q_i}_Q(M)\iso H^{q_i}_Q(M_i)\iso H^{q_i}_Q(M_i/M_{i-1})$,  where
$q_i=\cd(Q, M_i)$ for $i=1, \dots, r$ and $H^{k}_Q(M)=0$ for $k\not
\in \{ q_1, \dots, q_r\}.$
 \end{Lemma}
\begin{proof}
We proceed by induction on the length $r$ of $\mathcal{F}$.
The case $r=1$  is obvious.  Now suppose $r\geq 2$ and that the
statement holds for sequentially Cohen--Macaulay  modules with
respect to $Q$ with filtrations of length less than $r$.  We want to
prove it for $M$ which is sequentially Cohen--Macaulay with respect
to $Q$ and has the Cohen--Macaulay filtration
$\mathcal{F}$ of length $r$. Notice
that $M_{r-1}$ which appears in the filtration $\mathcal{F}$ of $M$
is also sequentially Cohen--Macaulay with respect to $Q$. Thus by
the induction hypothesis we have $H^{q_i}_Q(M_{r-1})\iso
H^{q_i}_Q(M_i)\iso H^{q_i}_Q(M_i/M_{i-1})$ for $i=1, \dots, r-1$ and
$H^{k}_Q(M_{r-1})=0$ for $k\not \in \{ q_1, \dots, q_{r-1}\}.$ Now
the exact sequence $0\rightarrow M_{r-1} \rightarrow M \rightarrow
M/M_{r-1} \rightarrow 0$ yields $H^{q_r}_Q(M)\iso
H^{q_r}_Q(M_r/M_{r-1})$ and $H^{t}_Q(M)\iso H^{t}_Q(M_{r-1})$ for
$0\leq t < q_r.$ Therefore,  the desired result follows.
\end{proof}
\begin{Fact}{\em
\label{grade}
In the proof of Lemma \ref{reg} one observes that
\[
\grade(Q, M_i)=q_1 \quad \text{for}\quad   i=1, \dots, r.
\]}
\end{Fact}
\begin{Theorem}
\label{joint} Let $L$ and $N$ be two non-zero finitely generated
graded modules over $K[x]$ and $K[y]$, respectively.  We set $M=L\tensor_KN$. Then the following statements are equivalent:
 \begin{itemize}
\item[{(a)}]  $M$ is a sequentially Cohen--Macaulay $S$-module with  respect to $Q$;

\item[{(b)}] $N$ is a sequentially Cohen--Macaulay $K[y]$-module.
\end{itemize}
\end{Theorem}
\begin{proof}
$(a)\Rightarrow (b)$: Let $\mathcal{F}$:
 $0=M_0\varsubsetneq M_1 \varsubsetneq  \dots \varsubsetneq M_r=M$ be the Cohen--Macaulay filtration with respect to $Q$ . By Lemma \ref{reg} we have
 \[
 H^{q_i}_Q(M)\iso H^{q_i}_Q(M_i)\iso H^{q_i}_Q(M_i/M_{i-1})
 \]
  where
 $q_i=\cd(Q, M_i)=\cd(Q, M_i/M_{i-1})$ for $i=1, \dots, r$ and $H^{k}_Q(M)=0$ for $k\not
\in \{ q_1, \dots, q_r\}.$
Note that  $
H^{q_i}_Q(M)\iso L\tensor_KH^{q_i}_Q(N)
$ for $i=1, \dots, r$, see also the proof of \cite[Proposition 1.5]{AR2}. Hence  $H^{q_i}_Q(M)^\vee \iso L^\vee \tensor_KH^{q_i}_Q(N)^\vee$  where $(-)^\vee$ is the Matlis-dual, see \cite[Lemma 1.1]{HR}. We conclude that
\begin{eqnarray*}
\big(H^{q_i}_Q(M_i/M_{i-1})^\vee\big)_{(k, *)} &\iso&  \big( H^{q_i}_Q(M)^\vee\big)_{(k, *)}\\
                                                              &\iso &(L^\vee)_k\tensor_KH^{q_i}_Q(N)^\vee \\
                                                              &\iso & \Ext^{n-q_i}_{K[y]}(N, K[y])^t,
\end{eqnarray*}
where $t=\dim_K(L^\vee)_k.$
Since each $M_i/M_{i-1}$ is Cohen--Macaulay with respect to $Q$ with $\cd(Q, M_i/M_{i-1})=q_i$, it follows from the above isomorphisms and Lemma \ref{cmd}   that  $\Ext^{n-q_i}_{K[y]}(N, K[y])$ is Cohen--Macaulay of dimension $q_i$ for $i=1, \dots, r$. If $k\not\in \{ q_1, \dots, q_r\}$, then  $ L\tensor_KH^{k}_Q(N) \iso H^{k}_Q(M)=0$ and hence $H^{k}_Q(N)=0$. Thus  $\Ext^{n-k}_{K[y]}(N, K[y])=0$  for $k\not\in \{ q_1, \dots, q_r\}$.   Therefore, the result follows from \cite[Theorem 1.4]{HS}.

  $(b)\Rightarrow (a)$: Let $N$ be  sequentially
Cohen--Macaulay $K[y]$-module with the Cohen--Macaulay filtration
$0=N_0\varsubsetneq N_1 \varsubsetneq \dots \varsubsetneq N_r=N$.
Consider the filtration $0=L\tensor_KN_0\subseteq L\tensor_KN_1
\subseteq \dots \subseteq L\tensor_KN_r=L\tensor_KN$. We claim this
filtration is the Cohen--Macaulay filtration with respect to
$Q$. First, we note that $L\tensor_KN_i  \varsubsetneq
L\tensor_KN_{i+1}$ for all $i$. Otherwise,   we have $\dim N_i=\dim
N_{i+1}$ by \cite[Corollary 2.3]{STY}, a contradiction.  For all $k$
and $i$  we have the following isomorphisms
\begin{eqnarray*}
H^k_Q\big( (L\tensor_KN_i)/(L\tensor_KN_{i-1}) \big)&\iso& H^k_Q(
L\tensor_K(N_i/N_{i-1}) \big )\\ &\iso& L\tensor_K
H^k_Q(N_i/N_{i-1}).
\end{eqnarray*}
The first isomorphism is standard,  and for the second one,  see the proof of \cite[Proposition 1.5]{AR2}.  We set $D_i= (L\tensor_KN_i)/(L\tensor_KN_{i-1})$ for all $i$. Thus
we have  $\cd (Q, D_i)=\dim N_i/N_{i-1}$ for all $i$. This implies
that  $\cd\big (Q, D_{i-1})<\cd\big (Q, D_i)$ for all $i$. Also each
 $D_i$ is Cohen--Macaulay with respect to $Q$
because $N_i/N_{i-1}$ is Cohen--Macaulay for all $i$.
\end{proof}

\section{Hypersurface rings which are sequentially Cohen--Macaulay with respect
to $Q$ } Let $f\in S$ be a bihomogeneous element of degree $(a, b)$
and consider the hypersurface ring $R=S/fS$. We may write
\[
f=\sum_{{| \alpha|=a}\atop {| \beta|=b}}c_{\alpha \beta }x^\alpha y^\beta  \;\; where \;\; c_{\alpha \beta} \in K.
\]
Notice that $R$ is a Cohen--Macaulay module of dimension $m+n-1.$  We collect some observations in the following
\begin{Lemma}
\label{hyper}
Consider the hypersurface ring $R$ defined above. Then the following statements hold:
\begin{itemize}
\item[{(a)}] If $a=0$ and $ b>0$, then $R$ is Cohen--Macaulay with respect to $P$ of $\cd(P, R)=m$ and Cohen--Macaulay with respect to $Q$ of $\cd(Q, R)=n-1.$
\item[{(b)}] If $a>0$ and $ b=0$, then $R$ is Cohen--Macaulay with respect to $P$ of $\cd(P, R)=m-1$ and Cohen--Macaulay with respect to $Q$ of $\cd(Q, R)=n.$
\item[{(c)}] If $ a>0$ and $b>0$, then $\grade(P, R)=m-1$ and $\cd(P, R)=m$, and $\grade(Q, R)=n-1$ and $\cd(Q, R)=n$.
\end{itemize}
\end{Lemma}
\begin{proof}
In order to proof (a), if $a=0$, then we may write $f=\sum_{| \beta|=b} c_\beta y^\beta$. Fact \ref{cd}(a) implies  $\cd(P, R)=\dim S/(Q+ (f))=m$ and $\cd(Q, R)=\dim S/(P+(f))=n-1.$
On the other hand, by Fact \ref{cd}(b), we have $\grade(P, R)= \dim R-\cd(Q, R)=m+n-1-(n-1)=m$ and $\grade(Q, R)= \dim R-\cd(P, R)=m+n-1-m=n-1$. Thus the conclusions follows. Parts (b) and (c) are proved in the same way.
\end{proof}
Notice that if $a, b>0$, then  $R$ is not Cohen--Macaulay
with respect to $Q$. Thus it is natural to ask whether $R$ is
sequentially Cohen--Macaulay with respect to $Q$.  In the following,
we classify all hypersurface rings that are sequentially
Cohen--Macaulay with respect to $Q$. We first show the following
\begin{Proposition}
\label{hyp3} Let $f\in S$ be a bihomogeneous element of degree $(a,
b)$ such that $f=h_1h_2$ where $h_1= \sum_
{\left|\alpha\right|=a}c_{\alpha} x^\alpha $ with $ c_{\alpha}\in K$
and $h_2= \sum_ {\left|\beta\right|=b}c_{\beta} y^\beta $ with $
c_{\beta}\in K$, i.e.,  $\deg h_1=(a, 0)$ and $\deg
h_2=(0,b)$. Consider the hypersurface ring $R=S/fS$. Then $R$ is
sequentially Cohen--Macaulay with respect to $P$ and $Q$.
\end{Proposition}
\begin{proof}
 We show that $R$ is sequentially Cohen--Macaulay with respect to $P$. The argument for $Q$ is similar.  Consider the filtration $\mathcal{F}$: $0=R_0\varsubsetneq R_1
\varsubsetneq R_2=R$ where $R_1=h_2S/fS$. We claim that this
filtration is the  Cohen--Macaulay filtration with respect
to $P$. Observe that $R_2/R_{1}\iso S/h_2S$ is
Cohen--Macaulay with respect to $P$ with $\cd(P, R_2/R_{1})=m$, by Lemma \ref{hyper}(a) . Now
consider the map $\varphi$: $S\longrightarrow h_2S/fS$ given by  $g\longmapsto gh_2+fS$. We get the
following isomorphism $S/h_1S\iso h_2S/fS\iso R_1/R_0$. Thus
$R_1/R_0$ is Cohen--Macaulay with respect to $P$ with
$\cd(P, R_1/R_0)=m-1$, by Lemma \ref{hyper}(b).  Therefore, $\mathcal{F}$ is the
Cohen--Macaulay filtration of $R$ with respect to $P$.
\end{proof}
For the proof of the main theorem, we recall the following results from \cite{PR}.
\begin{Fact}{\em
\label{PR}
 Let $\mathcal{D}: 0 = D_0\varsubsetneq   D_1\varsubsetneq\ldots\varsubsetneq D_r = M$ be the dimension filtration of $M$ with respect to $Q$. Then
 \begin{itemize}
\item[{(a)}]  	$D_i=\bigcap_{\pp_j\not \in B_{i,Q}}N_j$
	for $ i=1, \ldots, r-1$ where $0 =\bigcap_{j=1}^sN_j $ is a reduced primary decomposition of $0$
	in $M$ with $N_j$ is $\pp_j$-primary for $j = 1, \ldots, s$ and
	\[
B_{i,Q}=\{ \pp\in \Ass(M): \cd(Q, S/\pp)\leq \cd(Q, D_i)\}.
\]
\item[{(b)}]	$\Ass(M/D_i)=\Ass(M)\setminus \Ass(D_i)$  for $i=1,\ldots,r$.
\item[{(c)}]
$\grade(Q,M/D_{i-1}) = \cd(Q, D_i)$ for $i = 1, \ldots , r$  if and only if	$M$ is sequentially Cohen--Macaulay with respect to $Q$.
\end{itemize}}
\end{Fact}

\begin{Theorem}
\label{hypersurface} Let $f\in S$ be a bihomogeneous  element of
degree $(a, b)$ and  $R=S/fS$ be the hypersurface ring. Then the
following statements are equivalent:
 \begin{itemize}
\item[{(a)}] $R$ is sequentially Cohen--Macaulay with respect to
$Q$;
\item[{(b)}] $f=h_1h_2$ where $\deg h_1=(a, 0)$ with $a\geq 0$
and $\deg h_2=(0,b)$ with $b\geq 0$.
\end{itemize}
\end{Theorem}
\begin{proof}
 $(a)\Rightarrow (b)$:  We may assume that $R$ is not Cohen--Macaulay with respect to $Q$, see Lemma \ref{hyper}. Let $f=\prod_{i=1}^r f_i$ be the unique
factorization of $f$ into bihomogeneous irreducible factors  $f_i$
with $\deg f_i=(a_i, b_i)$ for $i=1, \dots, r$. Note that
$\sum_{i=1}^r a_i=a$
 and $\sum_{i=1}^r b_i=b$. Our aim is to show that for each $f_i$ we have $\deg f_i=(a_i,
 0)$ with $a_i\geq 0$ or $\deg f_i=(0 ,b_i)$  with $b_i\geq 0$. Assume this is not the case and
 so there exists $1\leq s \leq r$ such
 that $\deg f_s=(a_s ,b_s)$ with $a_s, b_s>0$.  Thus we may write
 that $\deg f_i=(a_i, 0)$ with $a_i\geq 0$ for $i=1, \dots, s-1$, $\deg f_i=(a_i ,b_i)$ with $a_i, b_i>0$  for $i=s,s+1, \dots, t$ and
 $\deg f_i=(0 ,b_i)$  with $b_i\geq 0$ for $i=t+1, \dots, r$, and $t<r$. By Fact \ref{PR}(a), $R$ has the
dimension filtration
$\mathcal{F}$:
$0=(f)/(f)\varsubsetneq I/(f) \varsubsetneq
  R=S/(f)$
with respect to $Q$  where  $I=\bigcap_{i=1}^t (f_i)$. Notice that $\cd(Q, I/(f))=n-1$ by Fact \ref{PR}(b) and Fact \ref{cd}(d), and $\cd(Q, R)=n$.
As $R$ is sequentially Cohen--Macaulay with respect to $Q$,  we have to have $\grade(Q, S/I)=\cd(Q,R)$ by Fact \ref{PR}(c). Since $S/I$ is Cohen--Macaulay, it follows from Fact \ref{cd}(b) that $\grade(Q, S/I)=\dim S/I-\cd(P, S/I)=(m+n-1)-m=n-1$, a contradiction.

 $(b) \Rightarrow (a)$:  follows from Proposition \ref{hyp3}.
\end{proof}

\begin{center}
{Acknowledgment}
\end{center}
\hspace*{\parindent}
The author would like to thank J\"urgen Herzog for his helpful comments.

\end{document}